\documentclass[11pt]{amsart}
\usepackage{amssymb,amsmath,latexsym,enumerate,graphicx}

\newtheorem*{rem}{Remark}

\newtheorem*{exam}{Example}

\def\Z{\mathbb{Z}}

\def\R{\mathbb{R}}

\def\H{\mathcal{H}}
\def\P{\mathcal{P}}

\title{Enumeration of $4 \times 4$ Magic Squares}

\author{Matthias Beck}
\address{Department of Mathematics\\
         San Francisco State University\\
         San Francisco, CA 94132\\
         U.S.A.}
\email{beck@math.sfsu.edu}

\author{Andrew Van Herick}
\address{531 Beloit Ave.\\
         Kensington, CA 94708\\
         U.S.A.}
\email{avanherick@gmail.com}

\thanks{
We thank a referee and an associate editor for helpful comments on an earlier version of this paper. 
We are grateful to San Francisco State University's Center for Computing and Life Sciences for graciously offering the use of their resources.
This research was partially supported by the NSF (research grant DMS-0810105).
}



\date{29 August 2009}

\begin{document}

\begin{abstract}
A \emph{magic square} is an $n \times n$ array of distinct positive integers whose sum along any row, column, or main diagonal is the same number. We compute the number of such squares for $n=4$, as a function of either the magic sum or an upper bound on the entries. The previous record for both functions was the $n=3$ case. Our methods are based on inside-out polytopes, i.e., the combination of hyperplane arrangements and Ehrhart's theory of lattice-point enumeration.
\end{abstract}

\maketitle


\section{It's A Kind Of Magic}

A \emph{magic square} is an $n \times n$ array of distinct positive integers whose sum along any row, column, or main diagonal is the same number, the \emph{magic sum}. The history of magic squares is well documented, see, e.g., \cite{CammannO,CammannI,zaslavskyclaudia}.
The contents of a magic square have varied with time and writer; usually they have been the first $n^2$ consecutive positive integers, but often any arithmetic sequence and sometimes fairly arbitrary numbers.  The fixed ideas are that they are integers, positive, and distinct.

In the last century mathematicians took an interest in results about the number of squares with a fixed magic sum, but with simplifications: diagonal sums were often omitted and the fundamental requirement of distinctness was almost invariably neglected \cite{ahmeddeloerahemmecke,magic,ehrhartmagic,stanleymagic}.
For example, classical formulas of MacMahon \cite{macmahon} include
\[
  3 {{ t+3 }\choose{ 4 }} + {{ t+2 }\choose{ 2 }} \ ,
\]
the number of $3\times 3$ squares with (not necessarily distinct) nonnegative integer entries that sum to $t$ along any row and column, and
\[
  \frac{2}{9} t^2 + \frac{2}{3} t + 1 \qquad \text{ if } 3|t ,
\]
the number of such squares in which the two main diagonals also sum to $t$.
The papers \cite{magiclatin,xinmagic} form, to the best of our knowledge, the beginning of a theory that tackles counting problems related to magic squares with the distinctness of the entries enforced.

Our goal is to show that the ideas in \cite{magiclatin} can be used to compute the number of $4 \times 4$ magic squares (with distinct entries), as a function of a parameter that is either the magic sum or an upper bound on the entries of the square.
To be precise, we define the \emph{affine magic counting function} $a_n(t)$ to be the number of all $n \times n$ matrices consisting of distinct positive integers whose sum along any row, column, or main diagonal is the same number $t$.
The \emph{cubical magic counting function} $c_n(t)$ is the number of $n \times n$ matrices whose entries are distinct positive integers less than $t$, whose sums along any row, column, or main diagonal are equal.
The previous record consisted of the $3 \times 3$ counting functions \cite{magiclatin,xinmagic} (see also \cite{sls} for the computational implementation of \cite{magiclatin})
\[
a_3(t) = \begin{cases}
\frac{1}{9} \left( 2t^2-32t+144 \right) &\text{if } t \equiv 0 \bmod 18 , \\
\frac{1}{9} \left( 2t^2-32t+78 \right) &\text{if } t \equiv 3 \bmod 18 , \\
\frac{1}{9} \left( 2t^2-32t+120 \right) &\text{if } t \equiv 6 \bmod 18 , \\
\frac{1}{9} \left( 2t^2-32t+126 \right) &\text{if } t \equiv 9 \bmod 18 , \\
\frac{1}{9} \left( 2t^2-32t+96 \right) &\text{if } t \equiv 12 \bmod 18 , \\
\frac{1}{9} \left( 2t^2-32t+102 \right) &\text{if } t \equiv 15 \bmod 18 , \\
0  &\text{if } t \not\equiv 0 \bmod 3
\end{cases}
\]
and~\cite{sls,magiclatin}
\[
c_3(t) = \begin{cases}
\frac{1}{6} \left( t^3-16t^2+76t-96 \right) &\text{if } t \equiv 0,2,6,8 \bmod 12 , \\
\frac{1}{6} \left( t^3-16t^2+73t-58 \right) &\text{if } t \equiv 1 \bmod 12 , \\
\frac{1}{6} \left( t^3-16t^2+73t-102 \right) &\text{if } t \equiv 3,11 \bmod 12 , \\
\frac{1}{6} \left( t^3-16t^2+76t-112 \right) &\text{if } t \equiv 4,10 \bmod 12 , \\
\frac{1}{6} \left( t^3-16t^2+73t-90 \right) &\text{if } t \equiv 5,9 \bmod 12 , \\
\frac{1}{6} \left( t^3-16t^2+73t-70 \right) &\text{if } t \equiv 7 \bmod 12 .
\end{cases}
\]
Both of these functions are \emph{quasipolynomials}, i.e., of the form $c_d(t) \, t^d + c_{ d-1 } (t) \, t^{ d-1 } + \dots + c_0(t)$, where $c_0, c_1, \dots, c_d$ are periodic functions. It follows from \cite{magiclatin} that $a_n(t)$ and $c_n(t)$ are always quasipolynomials in $t$; we will outline the basic arguments in Section \ref{iosection}.

A handy, compact way of representing a quasipolynomial $q(t)$ is through its generating function $Q(z) := \sum_{ t \ge 0 } q(t) \, z^t$. It is not hard to prove (see, e.g., \cite{ccd,stanleyec1}) that $Q(z)$ is a rational function with poles at $p$th roots of unity, where $p$ is a common period of the coefficient functions of $q(t)$. As examples we give the rational generating functions for $a_3(t)$ and $c_3(t)$:
\[
  A_3(z) = \frac{ 8 x^{15} \left( 2x^3+1 \right) }{ \left( 1-x^3 \right) \left( 1-x^6 \right) \left( 1-x^9 \right) }
\]
and
\[
  C_3(z) = \frac{ 8 x^{10} \left( 2 x^2 + 1 \right) }{ \left( 1-x \right)^2 \left( 1-x^4 \right) \left( 1-x^6 \right) } \, .
\]

Our computational approach (which differs from both \cite{magiclatin} and \cite{xinmagic}) verified the results for $a_3(t)$ and $c_3(t)$.
Both papers \cite{magiclatin,xinmagic} commented that the $4 \times 4$ case seems computationally infeasible, and we present some reasons for this assessment in Section \ref{compsection}; one of the reasons is that the rational generating functions for $a_4(t)$ and $c_4(t)$ take several pages to write down (in reduced form).
Nevertheless we were able to compute $A_4(z)$ and $C_4(z)$; since we did not want to waste paper, the results are posted online at
{\tt math.sfsu.edu/beck/papers/affmagic4.html} and {\tt math.sfsu.edu/beck/papers/cubmagic4.html}.


\section{Enter Geometry}\label{iosection}

We start with the affine case.
One treats a square with magic sum $t$ as an integer vector $x \in \Z^{n^2}$ confined to the affine subspace $ts$, where
\[
  s = \left\{ x \in \R^{n^2} : \text{ all line sums equal } 1 \right\} .
\]
From a geometric point, a square with magic sum $t$ is an integer point in the open dilated polytope $t\P = \R_{ >0 } \cap ts$, except that we need to require that the entries are distinct; so our integer point has to lie outside of the hyperplane arrangement
\[
  \H = \left\{ x_{ ij } = x_{ kl } : \, 1 \le i, j, k, l \le n, \, (i,j) \ne (k,l) \right\} .
\]
The pair $(\P, \H)$ is referred to as an \emph{inside-out polytope} \cite{iop}.
One can now use M\"obius inversion on this hyperplane arrangement (the computational approach in \cite{magiclatin}; one appealing feature of this approach is that one computes generating functions for polytopes of progressively smaller dimension) or view $\left( \R_{ >0 } \cap ts \right) \setminus \bigcup \H$ as a disjoint union of open polytopes (our approach, which we call \emph{region enumeration}). Either way, Ehrhart's fundamental results on integer-point enumeration in polytopes \cite{ccd,ehrhartpolynomial} imply that $a_n(t)$ is a sum/difference of \emph{Ehrhart quasipolynomials}---lattice-point enumerators of rational polytopes.

The cubical case is very similar, except that we replace the subspace by
\[
  s = \left\{ x \in \R^{n^2} : \text{ all line sums equal } \right\}
\]
and the polytope by $(0,1)^{ n^2 } \cap s$. The hyperplane arrangement, which enforces distinct entries in the square, is the same as in the affine case.


\section{Computations}\label{compsection}

For the affine case ($\P = \R_{ >0 } \cap s$), the present computational approach treats $t\P \setminus \bigcup \H$ as a disjoint union of open polytopes (the \emph{regions} of the inside-out polytope $(\P, \H)$) and computes $A_4(z)$ by summing the Ehrhart generating functions of all of the polytopes.  The computation of $C_4(z)$ in the cubical case is similar.  The challenge in this approach is computing the regions and summing the resulting rational functions efficiently.  Our algorithm for computing the regions of an inside-out polytope is fully described in \cite{vanherick} and runs in polynomial time for polytopes of fixed dimension.

To briefly outline the process, we observe that the removal of a single transverse hyperplane (i.e., one that intersects the relative interior of $\P$) from an open polytope \emph{splits} the polytope into two new open polytopes. Removing a non-transverse hyperplane simply results in the same open polytope. To enumerate the regions of an inside-out polytope, we begin by splitting the original polytope with the first hyperplane.  We then proceed recursively, splitting the newly generated polytope(s) with the next hyperplane, and so on, until all the hyperplanes in the arrangement are exhausted.  To determine whether a hyperplane is transverse to a given polytope, we engage a linear solver, using a normal vector to the hyperplane as the objective function, and the polytope as the feasible region.  Computation of each region's Ehrhart generating function is accomplished using Barvinok's algorithm \cite{barvinokehrhart}.


Aside from the running time of Barvinok's algorithm (which is polynomial in fixed dimension), the computational complexity of our approach has a similar order of magnitude to  M\"obius inversion, namely $\sum_{i = 1}^{d+1}{{ n }\choose{ i }}$ where $d = \dim \P$ and $n = |\H|$.  Implemented as a depth-first traversal of a binary tree \cite{vanherick}, memory requirements are linearly proportional to $|\H|$.  Despite the disadvantage of needing to compute Ehrhart generating functions for maximal-dimensional polytopes (as opposed to the M\"obius inversion approach), region enumeration lends itself well to parallel computation.  At a predetermined level of recursion, we send a description of each open polytope, together with the remaining hyperplanes, to separate processors and accumulate the results.

We have implemented the region enumeration algorithm using C++, combining the polytope libraries {\tt barvinok} \cite{verdoolaege}, {\tt polylib} \cite{polylib}, and {\tt cdd} \cite{cdd} with the computer algebra system {\tt GP/Pari} \cite{pari} to form a suite of applications and a small library capable of operations on inside-out polytopes, collectively referred to as {\tt IOP} \cite{vanherick} ({\tt http://iop.sourceforge.net/}).  As far as we are aware, {\tt IOP} is the first software implementation specifically capable of computing generating functions of inside-out polytopes.  The library includes the ability to convert an inside-out polytope to an isomorphic, lattice-point equivalent problem with full dimension in the ambient space, a feature that significantly reduces the time required to generate the regions.

To compute the number of $4 \times 4$ magic squares, we have reduced the scope of each problem by taking advantage of the $32$ symmetries described in \cite{vanherick}.  The reduced problems were fed to a 40-node Dell PowerEdge cluster. Computing the generating functions $A_4(t)$ and $C_4(t)$ each took less than a day. Each reduced problem required the enumeration of $3211412$ polytopes and the computation, summation, and simplification of an equal number of rational generating functions.  In lowest terms, the rational generating functions $A_4(t)$ and $C_4(t)$ respectively occupy approximately 280k and 25k of disk space with numerator/denominator polynomials of degree $2900$ (affine) and $533$ (cubical) and coefficients exceeding~$10^{42}$.  Files containing the results are available at {\tt math.sfsu.edu/beck/papers/[affmagic4.html,cubmagic4.html]}. The extreme nature of the simplified results justifies the skepticism voiced in \cite{magiclatin} and \cite{xinmagic} about the feasibility of computing counting functions for $4 \times 4$ magic squares.

During development of {\tt IOP} we observed phenomena that may be worthy of future investigation.  The degree to which rational-function arithmetic affected computation time was surprising. In a test run for the affine case, a 1.4 gHz home computer successfully enumerated the regions for a portion of the reduced problem \cite[Equation 6.16]{vanherick} and output individual simplified Ehrhart generating functions to a file in approximately a week's continuous run-time.  However, the attempt to generate a simplified sum proved intractable, with a projected run-time of well over three months.  This is considerably more gHz-hours than were actually used by the Dell cluster to tackle the entire reduced problem, suggesting that associative grouping of the rational generating functions may play a significant role in the run-time.  In particular, we suspect that grouping rational functions in expressions corresponding to an entire inside-out polytope may result in smaller simplified forms (hence, faster computation time) than expressions grouped arbitrarily.


\bibliographystyle{amsplain}

\def\cprime{$'$} \def\cprime{$'$}
\providecommand{\bysame}{\leavevmode\hbox to3em{\hrulefill}\thinspace}
\providecommand{\MR}{\relax\ifhmode\unskip\space\fi MR }
\providecommand{\MRhref}[2]{%
  \href{http://www.ams.org/mathscinet-getitem?mr=#1}{#2}
}
\providecommand{\href}[2]{#2}

\setlength{\parskip}{0cm}

\end{document}